\documentclass[12pt,twosided,reqno]{amsart}
\usepackage{amsmath}
\usepackage{amsthm}
\usepackage{amsfonts}
\usepackage{amssymb}

\theoremstyle{theorem}

\theoremstyle{definition}

\begin{document}

\vskip .5cm
\title{On the Archimedean characterization of parabolas}
 \vskip0.5cm
\thanks{
    \newline\indent { $^1$  supported by Basic Science Research Program through the National Research Foundation of Korea (NRF) funded by the Ministry of Education, Science and Technology (2010-0022926).}
\newline\indent { $^2$   supported by Basic Science Research Program through the National Research Foundation of Korea (NRF)
funded by the Ministry of Education, Science and Technology (2012R1A1A2042298) and supported by Kyungpook
National University Research Fund, 2012.}}

\vskip 0.5cm

\maketitle

\vskip 0.5cm
\centerline{\scshape
Dong-Soo Kim$^1$ and Young Ho Kim$^2$} \vskip .2in

\begin{abstract}
Archimedes knew that the area between a parabola and any chord $AB$
on the parabola is four thirds of the area of triangle $\Delta ABP$ where P is the
point on the parabola at which the tangent is parallel to $AB$. We
consider whether this property (and similar ones) characterizes parabolas.
We present five conditions  which are necessary and sufficient for
a strictly convex curve in the plane to be a parabola.
\end{abstract}

\vskip 1cm

\date{}
\maketitle

\section{Introduction}

 \vskip 0.50cm

A parabola is the set of points in the plane which are equidistant
from a point $F$ called the focus and a line $l$ called the directrix.
 Archimedes found some interesting area properties of parabolas.

Consider the region bounded by a parabola and a chord $AB$.
Let $P$ be the point on the parabola where the tangent is parallel to the chord $AB$.
The line through $P$ parallel to the axis of the parabola meets chord $AB$ at a point
$V$.
Then,  he showed that  the area of the parabolic region is
$a{|PV|}^{3/2}$ for some constant $a$, which depends only on the parabola.

Furthermore, he proved that the area of the parabolic region is $4/3$ times
the area of triangle $\bigtriangleup ABP$ whose base is the chord and whose third vertex is $P$.
For the proofs of Archimedes, see Chapter 7 of [8].
\vskip 0.3cm

In this paper, we
consider whether this property (and similar ones) characterizes parabolas.
As a result, we present five conditions  which are necessary and sufficient for
a strictly convex curve in the plane to be a parabola.

 \vskip 0.3cm
Usually, a curve $X$ in the  plane
 ${\mathbb R}^{2}$ is called {\it convex} if it bounds a convex domain in the  plane
 ${\mathbb R}^{2}$.

Hereafter,  we will  say  that a convex curve $X$ in the  plane
 ${\mathbb R}^{2}$ is  {\it strictly convex} if the curve   is smooth (that is, of class $C^{(3)}$) and is of positive  curvature $\kappa$
 with respect to the unit normal $N$ pointing to the convex side.
 Hence, in this case we have $\kappa(s)=\left< X''(s), N(X(s))\right> >0$, where $X(s)$ is an arclength parametrization of $X$.

For a smooth function $f:I\rightarrow {\mathbb R}$ defined on an open interval, we will also say that $f$ is
{\it strictly convex} if the graph of $f$ has  positive curvature $\kappa$ with respect to the upward unit normal $N$. This condition is equivalent to the positivity of $f''(x)$ on $I$.
\vskip 0.3cm

First of all, we prove the following characterization of parabolas:
 \vskip 0.3cm

 \noindent {\bf Theorem 1.}
 Let $X$ be the graph of a strictly convex function  $f:I\rightarrow {\mathbb R}$ in the  plane
 ${\mathbb R}^{2}$. Then $f$ is  a quadratic polynomial  if and only if   $X$ satisfies Condition:

 \vskip 0.3cm
 \vskip 0.3cm
 \noindent {\rm(A):}
 For a point $P$ on $X$ and a chord $AB$ of $X$ parallel to the tangent of $X$ at $P$,
 let $V$ denote  the point where the line through $P$ parallel to the $y$-axis meets $AB$.
 Then the area of the region bounded by the curve and $AB$ is $a{|PV|}^{3/2}$, where
 $a$ is a positive constant which depends only on the curve $X$.
 \vskip 0.3cm

Second, we prove
 \vskip 0.3cm

 \noindent {\bf Theorem 2.}
 Let $X$ be the graph of a strictly convex function  $f:I\rightarrow {\mathbb R}$ in  the plane
 ${\mathbb R}^{2}$. Then $f$ is  a quadratic polynomial  if and only if   $X$ satisfies Condition:

 \vskip 0.3cm
 \noindent {\rm(B):}
 For a sufficiently small  $k>0$, let $X_k$ denote the graph of $y=f(x)+k$.
 For any  point $V$ on $X_k$, let the tangent at $V$ meet the curve $X$ at $A$ and $B$.
 Then the region $S$ bounded by $X$ and the chord $AB$  has constant area (say, $\phi(k)$)
 independent of the choice of $V$.

 \vskip 0.3cm

Since  $|PV|=k$,  Theorem 1 is a special case of Theorem 2 for $\phi(k)=ak^{3/2}$, where $a$ is a constant.

 \vskip 0.3cm
Now, for an arbitrary strictly convex curve $X$ in the plane
 ${\mathbb R}^{2}$ which is not necessarily the  graph of a function,  we consider the following condition:

 \vskip 0.3cm
 \noindent {\rm(C):}
 For a point $P$ on $X$ and a chord $AB$ of $X$ parallel to the tangent of $X$ at $P$,
 the area of the region bounded by the curve and $AB$ is $4/3$ times
the area of triangle $\bigtriangleup ABP$.
 \vskip 0.3cm

Then,  we prove the following characterization of parabolas,
which is the  main theorem of this article.
\vskip 0.3cm

 \noindent {\bf Theorem 3.}
 Let $X$ be a strictly convex curve  in the plane
 ${\mathbb R}^{2}$. Then $X$ is a parabola if and only if  it satisfies Condition (C).

 \vskip 0.3cm

In order to prove Theorems 1, 2 and 3,  first of all,  in Section 2 we establish a new geometric meaning of  curvature $\kappa$ of a plane convex curve  $X$ at a point $P\in   M$ with $\kappa(P) >0$ (Lemma 6).
For the curvature function $\kappa$ of a plane curve, we refer to [3].

\vskip 0.3cm
As applications of  Theorem 3, we may prove some generalizations of Theorems 1 and 3 as follows.
\vskip 0.3cm

 \noindent {\bf Corollary 4.}
 Let $X$ be a strictly convex curve  in the plane
 ${\mathbb R}^{2}$. Then $X$ is a parabola if and only if  it satisfies Condition:

 \vskip 0.3cm
 \noindent {\rm(D):}
 For a point $P$ on $X$ and a chord $AB$ of $X$ parallel to the tangent of $X$ at $P$,
 the area of the region bounded by the curve and $AB$ is $a(P)|\bigtriangleup ABP|^{b(P)}$,
  where $a(P)$ and $b(P)$ are some functions of $P$ and $|\bigtriangleup ABP|$ denotes the area of the triangle
  $\bigtriangleup ABP$.
 \vskip 0.3cm

Finally, for the graph $X$ of a strictly convex function  $f:I\rightarrow {\mathbb R}$ in the plane
 ${\mathbb R}^{2}$, we consider the following Condition:
\vskip 0.3cm
 \noindent {\rm(E):}
 For a point $P$ on $X$ and a chord $AB$ of $X$ parallel to the tangent of $X$ at $P$,
 let $V$ denote  the point where the line through $P$ parallel to the $y$-axis meets $AB$.
 Then the area of the region bounded by the curve and $AB$ is $a(P) {|PV|}^{b(P)}$, where
 $a(P)$ and $b(P)$ are some functions of $P$.
\vskip 0.3cm
Then we prove
\vskip 0.3cm
 \noindent {\bf Corollary 5.}
 Let $X$ be the graph of a strictly convex function  $f:I\rightarrow {\mathbb R}$ in the plane
 ${\mathbb R}^{2}$. Then  $X$ satisfies Condition (E) if and only if  $X$ is a parabola, which is given by either
 a quadratic polynomial $f$ or a function $f$ in (3.26) according as the function $a(P)$ is constant or not.

 \vskip 0.3cm
It follows from Corollary 5 that Theorem 1 is a corollary of Theorem 3.

To prove Corollaries 4 and 5, first of all, we show that $b(P)$ must be $1$
 in  Corollary 4 (respectively, $3/2$ in  Corollary 5).
 Then we can show that $X$ satisfies Condition (C).
  Hence,   it follows from Theorem 3 that Corollaries 4 and 5 hold.
 \vskip 0.3cm

Among  the  graphs of functions,  \'A. B\'enyi et al. proved some characterizations of parabolas ([1,2])
and  B.  Richmond and T. Richmond  established a dozen characterizations of parabolas
using elementary techniques ([7]). In their papers, parabola means the graph of a quadratic polynomial in one
variable.

For an example, consider   a function $f(x)=b\{(1-cx)-\sqrt{1-2cx}\}$
in (3.26) with $b,c>0$ defined on $I=(-\infty, \frac{1}{2c})$.
Then, the function $f$ is strictly convex and its graph $X$ satisfies Condition (C) (but neither (A) nor (B)).
Note that $X$ is not
the graph of a quadratic polynomial, but an open part of the parabola given in (3.27).

\vskip 0.3cm
  Throughout this article, all curves are of class $C^{(3)}$ and connected, unless otherwise mentioned.
  \vskip 0.50cm

 \section{Preliminaries and  Theorem 1 and 2}
 \vskip 0.50cm

 Suppose that $X$ is a strictly convex  curve  in the plane
 ${\mathbb R}^{2}$ with the unit normal $N$ pointing to the convex side.
 For a fixed point $P \in X$, and for a sufficiently small $h>0$, consider the  line $l$ passing through
 $P+hN(P)$ which is parallel to
 the tangent of $X$ at $P$.
 Let's denote by $A$ and $B$ the points where the line $l$ intersects the curve $X$.

 We denote by $S_P(h)$ (respectively, $R_P(h)$) the area of the region bounded by the curve $X$ and chord $AB$
 (respectively, of the rectangle with a side $AB$ and another one on the tangent of $X$ at $P$ with height $h>0$).
 We also denote by $L_P(h)$ the length $|AB|$ of the chord $AB$.
  Then we have $R_P(h)=hL_P(h)=2|\bigtriangleup ABP|$, where $|\bigtriangleup ABP|$ denotes the area of the triangle $\bigtriangleup ABP$.

 We may adopt  a coordinate system $(x,y)$
 of  ${\mathbb R}^{2}$ in such a way that  $P$ is taken to be the origin $(0,0)$ and the $x$-axis is the tangent line of $X$ at $P$.
 Furthermore, we may assume that $X$ is locally  the graph of a non-negative strictly convex  function $f: {\mathbb R}\rightarrow {\mathbb R}$.
\vskip 0.3cm

For a sufficiently small $h>0$, we have

   \begin{equation}\tag{2.1}
   \begin{aligned}
   S_P(h)&=\int _{f(x)<h}\{h-f(x)\}dx,\\
   R_P(h)&= hL_P(h)=h\int _{f(x)<h}1dx.
       \end{aligned}
   \end{equation}
The  integration is taken on the interval $I_P(h)=\{x\in {\mathbb R}| f(x)<h\}$.

On the other hand,  we also have
 $$
   S_P(h)=\int _{y=0}^{h}L_P(y)dy.
     $$
This shows that
\begin{equation}\tag{2.2}
   \begin{aligned}
 S_P'(h)=L_P(h), \quad \text{and thus} \quad R_P(h)=hS_P'(h).
        \end{aligned}
   \end{equation}

 \vskip 0.3cm
First of all, we  prove the following lemma, which acts a key role in this article.
\vskip 0.3cm

 \noindent {\bf Lemma 6.} Suppose that   $X$  is a strictly convex  curve  in the plane
 ${\mathbb R}^{2}$ with  the unit normal $N$ pointing to the convex side.  Then
 we have
  \begin{equation}\tag{2.3}
   \begin{aligned}
   \lim_{h\rightarrow 0} \frac{1}{\sqrt{h}}L_P(h)= \frac{2\sqrt{2}}{\sqrt{\kappa(P)}},
    \end{aligned}
   \end{equation}
 where $\kappa(P)$ is the curvature of $X$ at $P$ with respect to  the unit normal $N$.
 \vskip 0.3cm

\noindent {\bf Proof.}
As above, we may adopt  a coordinate system $(x,y)$
 of  ${\mathbb R}^{2}$ in such a way that  $P$ is taken to be the origin $(0,0)$ and $X$ is locally  the graph of a non-negative strictly convex  function $f: {\mathbb R}\rightarrow {\mathbb R}$ with $f(0)=f'(0)=0$.
 Then $N$ is the upward unit normal.

  Since the curve $X$ is of class $C^{(3)}$, the Taylor's formula of $f(x)$ is given by
 \begin{equation}\tag{2.4}
 f(x)= ax^2 + f_3(x),
  \end{equation}
where  $a=f''(0)/2$, and $f_3(x)$ is an $O(|x|^3)$  function.
From  $\kappa(P)=f''(0)>0$, we see that $a$ is positive.

Now, we let $x=\sqrt{h}\xi$. Then, together with (2.1), (2.4) gives
\begin{equation}\tag{2.5}
   \begin{aligned}
  \frac{1}{\sqrt{h}}L_P(h)&= \frac{1}{\sqrt{h}}\int _{f(x)<h}1dx\\
  &=\int_{a\xi^2+  g_3(\sqrt{h}\xi)<1}1d\xi,
    \end{aligned}
   \end{equation}
where $g_3(\sqrt{h}\xi)=f_3(\sqrt{h}\xi)/h$.  Since $f_3$ is an  $O(|x|^3)$  function, we have
 \begin{equation}\tag{2.6}
   \begin{aligned}
   |g_3(\sqrt{h}\xi)|\le C\sqrt{h}|\xi|^3,
    \end{aligned}
   \end{equation}
where $C$ is a constant.
As $h \rightarrow 0$,
it follows from (2.5) and (2.6) that
\begin{equation} \tag{2.7}
  \begin{aligned}
\lim_{h\rightarrow 0} \frac{1}{\sqrt{h}}L_P(h)&=\int_{a\xi^2<1}1d\xi\\&= \frac{2}{\sqrt{a}}.
 \end{aligned}
  \end{equation}
Since $\kappa(P)=2a$, this completes the proof of Lemma 6.$\quad \square$

 \vskip 0.50cm
\noindent {\bf Remark.} From Lemma 6, we get  a new geometric meaning of  curvature $\kappa(P)$ of a plane convex curve  $X$ at a point $P\in   M$ with $\kappa(P) >0$. That is, we obtain
$$\kappa(P)=\lim_{h\rightarrow 0} \frac{8h}{L_P(h)^2}.$$

 \vskip 0.30cm
Now, we give a proof of Theorem 1.

 Let $X$ be  the graph of a strictly convex function  $f:I\rightarrow {\mathbb R}$, where $I$ is an open interval.
 Then $N$ is given by the upward unit normal.
For a fixed point
$P=(x,f(x))$ on $X$ and a small number $h>0$, consider the line $l$ passing through the point $P+hN(P)$
which is parallel to the tangent to $X$ at $P$.

Then the hypothesis shows that $S_P(h)=a|PV|^{3/2}$ for small $h>0$,
where $a$ is a constant depending only on $X$.
Note that  $|PV|=h\sec\theta$, where $f'(x)=\tan\theta$ is the slope of the tangent line at $P$.
Hence we have

\begin{equation} \tag{2.8}
  \begin{aligned}
S_P(h)&=a(\sec\theta)^{3/2} h^{3/2}\\
&=aW(x)^{3/2} h^{3/2},
  \end{aligned}
 \end{equation}
where $W(x)=\sqrt{1+f'(x)^2}$.
Thus  (2.2) yields
\begin{equation}\tag{2.9}
   \begin{aligned}
   L_P(h)=\frac{3}{2}aW(x)^{3/2} h^{1/2}.
    \end{aligned}
   \end{equation}
 \vskip 0.3cm
Therefore it follows from Lemma 6 that
\begin{equation}\tag{2.10}
   \begin{aligned}\kappa(P)=\frac{32}{9a^2W(x)^3}.
    \end{aligned}
   \end{equation}
Since the curvature $\kappa(P)$ of $X$ at $P=(x,f(x))$ is given by
\begin{equation}\tag{2.11}
   \begin{aligned}\kappa(P)=\frac{f''(x)}{W(x)^3},
    \end{aligned}
   \end{equation}
we see that $f''(x)$ is a constant.  Hence  $f(x)$ is a  quadratic polynomial.
This completes the proof of the if part of Theorem 1.

 By a straightforward calculation, it is trivial
to prove the only if part  of Theorem 1.
This completes the proof of Theorem 1.

 \vskip 0.50cm

Second,  we give a proof of Theorem 2.

Let $X$ be  the graph of a strictly convex function  $f:I\rightarrow {\mathbb R}$, where $I$ is an open interval.
 Then $N$ is given by the upward unit normal.
We fix a  point
$P(x,f(x))$ on $X$.  For a sufficiently small $h>0$,
consider the  line $l$ passing through
 $P+hN(P)$ which is parallel to
 the tangent of $X$ at $P$.
 Let's denote by $A$ and $B$ the points where the line $l$ intersects the curve $X$.

Then the chord $AB$ is tangent to $X_k$ at $V(x, f(x)+k)$, where $k=hW$ and  $W(x)=\sqrt{1+f'(x)^2}$.
The hypothesis shows that $S_P(h)=\phi(k)$.
It follows  from (2.2) that

\begin{equation} \tag{2.12}
  \begin{aligned}
L_P(h)&=S_P'(h)=W(x)\phi'(hW),\\
R_P(h)&=hL_P(h)=hW(x)\phi'(hW).\\
 \end{aligned}
 \end{equation}
Hence we have
\begin{equation} \tag{2.13}
  \begin{aligned}
\frac{L_P(h)}{\sqrt{h}}= \frac{\phi'(k)}{\sqrt{k}}W(x)^{3/2}.
 \end{aligned}
 \end{equation}

 For a fixed point
$P(x,f(x))$ on $X$, it follows from  $k=hW(x)$ that $h\rightarrow 0$ is equivalent to  $k\rightarrow 0$.
Thus,  Lemma 6 implies that
\begin{equation}\tag{2.14}
   \begin{aligned}
  \lim_{k\rightarrow 0}\frac{\phi'(k)}{\sqrt{k}}=
   W(x)^{-3/2}\lim_{h\rightarrow 0} \frac{1}{\sqrt{h}}L_P(h)= \frac{2\sqrt{2}}{\sqrt{\kappa(P)}}  W(x)^{-3/2}.
    \end{aligned}
   \end{equation}
If we denote by $\alpha$ the limit of the left hand side of (2.14), which is independent of $P$,  then we have
\begin{equation}\tag{2.15}
   \begin{aligned}\kappa(P)=\frac{8}{\alpha^2W(x)^3}.
    \end{aligned}
   \end{equation}
 \vskip 0.3cm

Similarly to the proof of Theorem 1,
we see that $f(x)$ is a quadratic polynomial.
This completes the proof of the if part of Theorem 2.

 For   a proof of the only if part of Theorem 2,  see Example 1.2 in [6, p.6].
This completes the proof of Theorem 2.

 \vskip 0.50cm

\section{Main Theorem}
 \vskip 0.3cm
In this section, we prove Theorem 3, which is the main theorem of this article.

  Let  $X$ denote a strictly convex  curve  in the plane
 ${\mathbb R}^{2}$ with the unit normal $N$ pointing to the convex side.
Suppose that  $X$ satisfies Condition (C). Then, for $P\in X$ and a sufficiently small $h>0$ we have
\begin{equation} \tag{3.1}
  \begin{aligned}
S_P(h)&=\frac{2}{3}R_P(h).
 \end{aligned}
 \end{equation}
 By differentiating (3.1) with respect to $h$, it follows from (2.2) that
\begin{equation}\tag{3.2}
   \begin{aligned}
   L_P(h)=2hL_P'(h).
    \end{aligned}
   \end{equation}
 \vskip 0.3cm
Therefore, we get
\begin{equation}\tag{3.3}
   \begin{aligned}
   L_P(h)=c(P)\sqrt{h},
    \end{aligned}
   \end{equation}
   where $c=c(P)$ is a constant depending on  $P$.
Furthermore,  Lemma 6 implies that
\begin{equation}\tag{3.4}
   \begin{aligned}
c(P)=\frac{2\sqrt{2}}{\sqrt{\kappa(P)}}.
 \end{aligned}
   \end{equation}
 \vskip 0.3cm
In order to prove Theorem 3, first, we  fix an arbitrary  point $A$ on $X$.

As before, we take a coordinate system $(x,y)$
 of  ${\mathbb R}^{2}$: $A$ is taken to be the origin $(0,0)$ and $x$-axis is the tangent line of $X$ at $A$.
 Furthermore, we may regard $X$ to be locally  the graph of a non-negative strictly convex  function $f: {\mathbb R}\rightarrow {\mathbb R}$
 with $f(0)=f'(0)=0$ and $f''(0)>0$.

For any point $B(x,f(x))$ with $x\ne 0$, we denote by $P$ the point on $X$ such that
the chord $AB$ is parallel to the tangent of $X$ at $P$. Then  we have
$P=(g(x), f(g(x)))$, for a function $g: {\mathbb R}\setminus \{0\}
 \rightarrow {\mathbb R}$ which satisfies $|g(x)|<|x|$ and

\begin{equation}\tag{3.5}
   \begin{aligned}
  xf'(g(x))=f(x).
 \end{aligned}
   \end{equation}
Since   $g(x)$ tends to $0$ as $x\rightarrow 0$, we may assume that $g(0)=0$.
 \vskip 0.3cm
We prove the following lemma, which plays a crucial role in the proof of Theorem 3.

\vskip 0.3cm

 \noindent {\bf Lemma 7.} $f(x)$ and $g(x)$ satisfy

   \begin{equation}\tag{3.6}
   \begin{aligned}
   x^3 f''(g(x))=8\{f(x)g(x)-xf(g(x))\},
    \end{aligned}
   \end{equation}

   \begin{equation}\tag{3.7}
   \begin{aligned}
   xf(x)=\frac{4}{3}\{f(x)g(x)-xf(g(x))\}+2\int_0^{x}f(t)dt.
    \end{aligned}
   \end{equation}
 \vskip 0.3cm

\noindent {\bf Proof.} Consider the triangle $\bigtriangleup ABC$, where $C$ denotes the point $(x,0)$. Then we have
$|AC|^2+|BC|^2=|AB|^2$. Note that  by definition, $|AB|^2=L_P(h)^2$, where $h$ denotes the distance from $P$ to the chord $AB$.
This shows that
 \begin{equation}\tag{3.8}
   \begin{aligned}
   x^2+f(x)^2=L_P(h)^2.
    \end{aligned}
   \end{equation}
The distance $h$ from $P$ to the chord $AB$ is given by
 \begin{equation}\tag{3.9}
   \begin{aligned}
  h=\frac{\epsilon \{f(x)g(x)-xf(g(x))\}}{\sqrt{x^2+f(x)^2}},
    \end{aligned}
   \end{equation}
where $\epsilon = 1$ for $x>0$ and $\epsilon = -1$ for $x<0$.
 \vskip 0.3cm
Since the curvature $\kappa(P)$ of $X$ at $P$ is given by
 \begin{equation}\tag{3.10}
   \begin{aligned}
  \kappa(P)=\frac{f''(g(x))}{(\sqrt{1+f'(g(x))^2})^3},
    \end{aligned}
   \end{equation}
it follows from (3.3), (3.4) and (3.5) that
 \begin{equation}\tag{3.11}
   \begin{aligned}
  L_p(h)^2=\frac{8h}{\kappa (P)}=\frac{8(x^2+f(x)^2)}{f''(g(x))x^3}\{f(x)g(x)-xf(g(x))\}.
    \end{aligned}
   \end{equation}
Together with  (3.8), this implies that (3.6) holds.

 \vskip 0.3cm
In order to prove (3.7), we consider the area of  triangle $\bigtriangleup ABC$. Then we have

 \begin{equation}\tag{3.12}
   \begin{aligned}
  \frac{\epsilon}{2}xf(x)=S_P(h)+\epsilon \int_0^{x}f(t)dt,
    \end{aligned}
   \end{equation}
   where $\epsilon = 1$ for $x>0$ and $\epsilon = -1$ for $x<0$.
By assumption, we have $S_P(h)=(4/3)|\bigtriangleup ABP|$. Hence we get
 \begin{equation}\tag{3.13}
   \begin{aligned}
 S_P(h)=\frac{2\epsilon}{3}\{f(x)g(x)-xf(g(x))\}.
    \end{aligned}
   \end{equation}
Together with (3.12), this implies that (3.7) holds. $\quad \square$
 \vskip 0.50cm

Next, with the help of Lemma 7, we show that in a neighborhood of an arbitrary point $A\in X$, the curve $X$ is a parabola.

By differentiating (3.7) with respect to $x$,
it follows from (3.5) that
 \begin{equation}\tag{3.14}
   \begin{aligned}
  f(g(x))=g(x)f'(x)-\frac{3}{4}\{xf'(x)-f(x)\}.
    \end{aligned}
   \end{equation}
Differentiating (3.5) with respect to $x$, and using again (3.5), we get
\begin{equation}\tag{3.15}
   \begin{aligned}
   f''(g(x))=\frac{xf'(x)-f(x)}{x^2g'(x)}.
    \end{aligned}
   \end{equation}

On the other hand, together with (3.14), (3.6) shows that
\begin{equation}\tag{3.16}
   \begin{aligned}
   f''(g(x))=\frac{xf'(x)-f(x)}{x^3}\{6x-8g(x)\}.
    \end{aligned}
   \end{equation}
It follows from (3.15) and (3.16) that
\begin{equation}\tag{3.17}
   \begin{aligned}
  \{xf'(x)-f(x)\}\{8g(x)g'(x)-6xg'(x)+x\}=0.
    \end{aligned}
   \end{equation}
Since  $f(x)$ is strictly convex, we obtain
\begin{equation}\tag{3.18}
   \begin{aligned}
  8 g(x)g'(x)-6xg'(x)+x=0.
    \end{aligned}
   \end{equation}

If we let $y=g(x)$, then (3.18) becomes $xdx+ (8y-6x)dy=0$.
 By putting $y=vx$, we get a separable differential equation, and hence
we can solve (3.18).
Since $g(0)=0$, we see that  $g(x)=x/2, x/4$ or
\begin{equation}\tag{3.19}
   \begin{aligned}
   g(x)=\frac{1}{4c}(cx+1-\sqrt{1-2cx}),
    \end{aligned}
   \end{equation}
where $c$ is a  nonzero constant.

By differentiating (3.14) with respect to $x$,
it follows from (3.5) that
 \begin{equation}\tag{3.20}
   \begin{aligned}
 \{xg(x)-\frac{3}{4}x^2\}f''(x)+xg'(x)f'(x)- g'(x)f(x)=0.
    \end{aligned}
   \end{equation}
 \vskip 0.3cm
If $g(x)=x/2$, then (3.20) shows that
  \begin{equation}\tag{3.21}
   \begin{aligned}
   x^2 f''(x)-2xf'(x)+2f(x)=0,
    \end{aligned}
   \end{equation}
of which  general solutions are given by $ax^2+bx$ for some $a,b\in {\mathbb R}$.
Since $f(0)=f'(0)=0$, it follows from (3.21) that $f(x)=ax^2$ for some positive constant $a$.
Thus, in a neighborhood of $A$, the curve $X$ is   a parabola.

 \vskip 0.3cm
If $g(x)=x/4$, then (3.20) yields that
  \begin{equation}\tag{3.22}
   \begin{aligned}
   2x^2 f''(x)-xf'(x)+f(x)=0.
    \end{aligned}
   \end{equation}
For some $a,b\in {\mathbb R}$, the general solutions of (3.22) are given by
\begin{equation}\tag{3.23}
   \begin{aligned}
   f(x)=ax+b\sqrt{|x|}.
    \end{aligned}
   \end{equation}
 This contradicts to   $f'(0)=0$.
\vskip 0.3cm

If $  g(x)=\frac{1}{4c}(cx+1-\sqrt{1-2cx})$, it follows from (3.20) that
 \begin{equation}\tag{3.24}
   \begin{aligned}
  (1-2cx)\{ \sqrt{1-2cx}-(1-cx)\}f''(x)+c^2xf'(x)-c^2f(x)=0.
    \end{aligned}
   \end{equation}
The general solutions of (3.24) are given by
 \begin{equation}\tag{3.25}
   \begin{aligned}
   f(x)=ax+b(1-\sqrt{1-2cx}),
    \end{aligned}
   \end{equation}
where  $a,b\in {\mathbb R}$.
Since $f(x)$ satisfies $f(0)=f'(0)=0$ and $f''(0)>0$, (3.25) shows that
 \begin{equation}\tag{3.26}
   \begin{aligned}
   f(x)=b\{(1-cx)-\sqrt{1-2cx}\},
    \end{aligned}
   \end{equation}
where $b$ is a positive constant. Hence, in a neighborhood of $A$, the curve $X$ is given by
 \begin{equation}\tag{3.27}
   \begin{aligned}
   b^2c^2x^2+2bcxy+y^2-2by=0.
    \end{aligned}
   \end{equation}

It follows from  the classification theorem of quadratic polynomials in $x$ and $y$ that the curve defined by (3.27)
 is a parabola.

Summarizing the above discussions, we see that the curve $X$ is locally a parabola.
 \vskip 0.3cm

 Finally, we show that the curve $X$ is a parabola as follows.

First, consider two parabolas $\Phi_1$  and $\Phi_2$ in the  plane
 ${\mathbb R}^{2}$.
For each $i=1,2$, let's denote by $\phi_i$  a connected open arc of the  parabola $\Phi_i$.

 Suppose that the two arcs $\phi_1$  and $\phi_2$ share a common subarc $\phi$.
 We fix a point $A$ on the subarc $\phi$.  As before, we take a coordinate system $(x,y)$
 of  ${\mathbb R}^{2}$: $A$ is taken to be the origin $(0,0)$, $x$-axis is the tangent line of $\phi$ at $A$
  and $\phi$ lies in the upper half plane. Then for each $i=1,2$,
 the parabolic arc $\phi_i$ is locally  the graph of $f_i$
 which is either of the form $f_i(x)=a_ix^2$ with $a_i>0$ or of the form in (3.26) with $b=b_i>0,c=c_i\ne0$.
  That is,  the parabola $\Phi_i$
  is of the form $y=a_ix^2$ with $a_i>0$ or of the form
 in (3.27) with $b=b_i>0,c=c_i\ne0$.

 Since $f_1$ is equal to $f_2$ around $x=0$, $f_1$ and $f_2$ have the same derivatives at the origin.
 Hence, we immediately see that
 $\Phi_1=\Phi_2$ because $f_i''(0)=2a_i, f_i'''(0)=0$ or $f_i''(0)=b_ic_i^2, f_i'''(0)=3b_ic_i^3$ in each case for $i=1,2$.

 Next, let's fix a point $A$ on the curve $X$. Then an open arc of $X$ containing $A$ is a  parabolic arc $\phi_0$ of a parabola $\Phi_0$. For an arbitrary point $B$ on the curve $X$,
   the compactness of the closed arc $AB$ of $X$ shows that  there exist  consecutive points $A=P_0, P_1, \cdots, P_n=B$ on $X$ and open  arcs $\phi_0, \phi_1, \cdots, \phi_n$ of $X$ such that
 1) for each $i=0,1,\cdots,n$, $P_i$ lies on $\phi_i$, 2) each $\phi_i$ is a parabolic arc,  3) $\{\phi_i\}$ covers the closed arc $AB$ of $X$.

 Since  $\phi_i$ and $\phi_{i+1}$ share a common subarc for  each $i=0,1,\cdots,n-1$,
 a successive use of the above argument shows that every $\phi_i$ is an arc of the parabola $\Phi_0$, and hence $B\in  \Phi_0$.
Therefore we see that $X$  is the parabola $\Phi_0$.

 This completes the proof of the if part of Theorem 3.

 For  a  proof of the only if part of Theorem 3,  see Chapter 7 of [8], which is originally due to Archimedes.
This completes the proof of Theorem 3.

 \vskip 0.50cm

\section{Corollaries and Remarks}
 \vskip 0.3cm

In this section, first of all, we prove Corollaries 4 and 5.

First, suppose that a strictly convex curve $X$ in the plane ${\mathbb R}^{2}$ satisfies Condition (D) with $b(P)=1$. Then we have
\begin{equation} \tag{4.1}
  \begin{aligned}
S_P(h)&=\frac{a(P)}{2}hL_P(h).
 \end{aligned}
 \end{equation}
 By differentiating (4.1) with respect to $h$, it follows from (2.2) that
\begin{equation}\tag{4.2}
   \begin{aligned}
   (2-a(P))L_P(h)=a(P)hL_P'(h).
    \end{aligned}
   \end{equation}
 \vskip 0.3cm
Solving (4.2), we get
\begin{equation}\tag{4.3}
   \begin{aligned}
   L_P(h)=c(P)h^{d(P)},
    \end{aligned}
   \end{equation}
 where    $c=c(P)$ is a constant depending on  $P$ and  $ d(P)=(2-a(P))/a(P)$.

It follows from (4.3) and  Lemma 6 that $d(P)=1/2$, and hence, $a(P)=4/3.$
Thus,  the curve $X$ satisfies Condition (C).

Now, suppose that  $X$ satisfies Condition (D) with $b(P)\ne 1$. Then we have
\begin{equation} \tag{4.4}
  \begin{aligned}
S_P(h)&=a(P)2^{-b(P)}\{hL_P(h)\}^{b(P)},
 \end{aligned}
 \end{equation}
 which shows that $b(P)>0$.
 By differentiating (4.4) with respect to $h$, it follows from (2.2) that
\begin{equation}\tag{4.5}
   \begin{aligned}
   L_P'(h)+h^{-1}L_P(h)=c(P)h^{-b(P)}L_P(h)^{2-b(P)},
    \end{aligned}
   \end{equation}
where $c(P)=2^{b(P)}a(P)^{-1}b(P)^{-1}$.
Solving the Bernoulli equation (4.5), we get
\begin{equation}\tag{4.6}
   \begin{aligned}
  \{h L_P(h)\}^{b(P)-1}=c(P)(b(P)-1)\ln h +d(P),
    \end{aligned}
   \end{equation}
   where $d(P)$ is a constant depending on  $P$.

In case $b(P)>1$, by letting $h\rightarrow0$, (4.6) leads to a contradiction.
In case $b(P)\in (0,1)$, multiplying the both sides of (4.6) by $h^{\alpha(P)}$ with $\alpha(P)=(1-b(P))/2>0$,
and then  by letting $h\rightarrow0$, we get a contradiction.
This  shows that  $b(P)$ must be $1$.

Together with the above discussion on the case $b(P)=1$, Theorem 3 completes the proof  of Corollary 4.

 \vskip 0.50cm

 Next, we prove Corollary 5.

  Suppose that  the graph  $X$ of a strictly convex function  $f:I\rightarrow {\mathbb R}$ in the plane
 ${\mathbb R}^{2}$ satisfies Condition (E). Then for a fixed point $P(x,f(x))$ on $X$ and for  $h>0$,
 we have
 \begin{equation}\tag{4.7}
   \begin{aligned}
   S_P(h)=a(P) {|PV|}^{b(P)}.
    \end{aligned}
   \end{equation}
 Since $|PV|=hW(x)$ with $W(x)=\sqrt{1+f'(x)^2}$, by differentiating (4.7) with respect to $h$, we get
  \begin{equation}\tag{4.8}
   \begin{aligned}
   L_P(h)=a(P)b(P) W(x)^{b(P)}h^{b(P)-1}.
    \end{aligned}
   \end{equation}

 Hence, it follows from Lemma 6 that $b(P)=3/2$. This shows that
  \begin{equation}\tag{4.9}
   \begin{aligned}
    S_P(h)=a(P) W(x)^{3/2}h^{3/2} \quad \text{and} \quad L_P(h)=\frac{3}{2}a(P) W(x)^{3/2}\sqrt{h}.
    \end{aligned}
   \end{equation}
Thus we get
 \begin{equation}\tag{4.10}
   \begin{aligned}
   \Delta ABP= \frac{1}{2}hL_P(h) =\frac{3}{4}a(P) W(x)^{3/2}h^{3/2}=\frac{3}{4} S_P(h),
    \end{aligned}
   \end{equation}
which shows that $X$ satisfies Condition (C). Therefore, it follows from the proof of Theorem 3
 that  $X$ is a parabola, which is given by either a quadratic polynomial $f$
or a function $f$ in (3.26).

Conversely, if $f$ is  a quadratic polynomial, Theorem 1 shows that the graph $X$ of $f$ satisfies
Condition (E) with  a constant $a(P)$ and $b(P)=3/2$.
If  $f$ is  a function  in (3.26), it is straightforward to show that the graph $X$ of $f$ satisfies
Condition (E) with  a nonconstant function $a(P)$ and $b(P)=3/2$.

This completes the proof of Corollary 5.

\vskip 0.3cm

Together with (3.1)-(3.4) and Theorem 3, the same argument  as in the proof of Corollary  5 shows
\vskip 0.3cm
\noindent {\bf Corollary 8.} Let $X$ denote   a strictly convex curve  in the plane ${\mathbb R}^{2}$.
Then, the following are equivalent.
\vskip 0.3cm
\noindent 1) $X$ satisfies Condition (C).

\noindent 2) $ S_P(h)=a(P)h^{3/2},$ where  $a(P)$ is a function of $P\in X$.

\noindent 3) $ S_P(h)=a(P)h^{b(P)},$ where  $a(P)$ and $b(P)$ are some  functions of $P\in X$.

\noindent 4) $X$ is a parabola.

\vskip 0.3cm
\noindent {\bf Remark 9.} It follows   from our proofs that Theorem 3 holds even if a strictly convex (hence, $C^{(3)}$) curve $X$  satisfies
Condition (C) for sufficiently small $h>0$
 at every point $P\in X$.
\vskip 0.3cm

 Finally, we give an example of a convex curve which satisfies Condition (C) for sufficiently small $h>0$ at every point $P\in X$,  but it is not a  parabola.
Note that the example is not of class $C^{(3)}$, and hence it is not strictly convex either.
\vskip 0.3cm

 \noindent {\bf Example 10.}
Consider the graph $X$ of a function $f:{\mathbb R}\rightarrow {\mathbb R}$ which is given by
\begin{equation}\tag{4.11}
 \begin{aligned}
 f(x)= \begin{cases}
  9x^2, & \text{if $x<0,$} \\
 \frac{9}{4}x^2, & \text{if $x\ge 0.$}
  \end{cases}
  \end{aligned}
  \end{equation}
 Then, the function $f$ is not of class  $C^{(3)}$ at the origin, and hence the curve $X$ is not strictly convex.
It is straightforward to show that if $P$ is the origin, then for all $h>0$ we have
\begin{equation}\tag{4.12}
   \begin{aligned}
  L_P(h)=\sqrt{h},  \quad \text{and} \quad S_P(h)=\frac{2}{3}R_P(h).
    \end{aligned}
   \end{equation}
Hence $X$ satisfies Condition (C) at the origin  for all $h>0$. If $P\in X$ is not the origin, then
there  exists a positive number $\varepsilon(P)$ such that
 for every positive number $h$ with $h<\varepsilon(P)$,
$X$ satisfies Condition (C).

Thus, $X$ satisfies Condition (C) for sufficiently small $h>0$ at every point $P\in X$. But it is not a parabola.
\vskip 0.3cm

\noindent {\bf Remark 11.}
In [4] and [5], the authors proved  the  higher dimensional versions of Theorems 1 and 2, respectively.
 \vskip 0.50cm


\vskip 1.0 cm

Department of Mathematics, \par Chonnam National University,\par
Kwangju 500-757, Korea

{\tt E-mail:dosokim@chonnam.ac.kr} \vskip 0.3 cm

Department of Mathematics, \par Kyungpook National University,\par
Taegu 702-701, Korea

{\tt E-mail:yhkim@knu.ac.kr} \vskip 0.3 cm


\vskip 0.50cm

\begin{thebibliography}{5.4}
\bibitem {}
B\'enyi, \'A., Szeptycki, P. and  Van Vleck, F., {\it Archimedean properties of parabolas}, Amer. Math. Monthly 107 (2000), no. 10, 945-949.

\bibitem {}
B\'enyi, \'A., Szeptycki, P. and  Van Vleck, F., {\it A generalized Archimedean property}, Real Anal. Exchange 29 (2003/04), no. 2, 881-889.
\bibitem {}
do Carmo, M. P., {\it Differential Geometry of
Curves and Surfaces}, Prentice-Hall, Englewood Cliffs,
NJ, 1976.

\bibitem {}
Kim, D.-S.  and  Kim, Y. H., {\it Some characterizations of spheres and elliptic paraboloids}, Linear Algebra Appl.
437 (2012), no. 1, 113-120.


\bibitem {}
Kim, D.-S.  and  Kim, Y. H., {\it Some characterizations of spheres and elliptic paraboloids II},
 Linear Algebra Appl., 438 (2013), no. 3, 1356-1364.

\bibitem {}
Nomizu, K. and Sasaki, T., {\it Affine differential geometry. Geometry of affine immersions.
 Cambridge Tracts in Mathematics, 111}, Cambridge University Press, Cambridge, 1994.

\bibitem {}
Richmond, B. and Richmond, T.,  {\it How to recognize a parabola}, Amer. Math. Monthly 116(2009), no.10, 910-922.

\bibitem {}
Stein, S., {\it Archimedes. What did he do besides cry Eureka?},
Mathematical Association of America, Washington, DC, 1999.

\end{thebibliography}
\end{document}